\documentclass[12pt]{article}
\usepackage{amsfonts}
\usepackage{bbm}
\usepackage{mathrsfs}
\usepackage{amssymb}

\title{On the Construction of Gr\"obner Bases with Coefficients\\  in Quotient Rings\thanks{Project supported by the National Natural
Science Foundation of China (10971044).}}

\author{Huishi Li\\
{\small Department of Applied Mathematics, College of Information Science and Technology}\\
{\small Hainan University, Haikou 570228, China}}

\date{}

\begin{document}
\maketitle
\begin{center}
\begin{minipage}{120mm}
{\small {\bf Abstract.} Let $\Lambda$ be a commutative Noetherian
ring, and let $I$ be a proper  ideal of $\Lambda$, $R=\Lambda /I$.
Consider the polynomial rings $T=\Lambda [x_1,\ldots x_n]$  and
$A=R[x_1,\ldots ,x_n]$. Suppose that linear equations are solvable
in $\Lambda$. It is shown that linear equations are solvable in $R$
(thereby theoretically Gr\"obner bases for ideals of $A$ are well
defined and constructible) and that practically Gr\"obner bases in
$A$  with respect to any given monomial ordering can be obtained by
constructing Gr\"obner bases in $T$, and moreover, all basic
applications of a Gr\"obner basis at the level of $A$ can be
realized by a Gr\"obner basis at the level of $T$. Typical
applications of this result are demonstrated respectively in the
cases where $\Lambda=D$ is a PID, $\Lambda =D[y_1,\ldots ,y_m]$ is a
polynomial ring over a PID  $D$, and  $\Lambda =K[y_1,\ldots ,y_m]$
is a polynomial ring over a field $K$. }
\end{minipage}\end{center} {\parindent=0pt\vskip 6pt

{\bf 2010 MSC} 13P10\vskip 6pt

{\bf Key words} Commutative ring, quotient ring,  Gr\"obner basis.}

\def\NZ{\mathbb{N}}
\def\QED{\hfill{$\Box$}}
\def \r{\rightarrow}
\def\normalbaselines{\baselineskip 24pt\lineskip 4pt\lineskiplimit 4pt}
\def\mapright#1#2{\smash{\mathop{\longrightarrow}\limits^{#1}_{#2}}}

\def\v5{\vskip .5truecm}
\def\T#1{\widetilde #1}
\def\H#1{\widehat #1}

\def\OV#1{\overline {#1}}
\def\hang{\hangindent\parindent}
\def\textindent#1{\indent\llap{#1\enspace}\ignorespaces}
\def\item{\par\hang\textindent}

\def\LM{{\bf LM}}\def\LT{{\bf
LT}}\def\B{{\cal B}} \def\LC{{\bf LC}} \def\G{{\cal G}} 

\vskip 1truecm

\section*{1. Introduction and Preliminary}
It is well known that the celebrated theory of Gr\"obner bases over 
a field invented by Buchberger ([Bu1], [Bu2]) has been generalized 
to study Gr\"obner bases and their applications over various 
commutative rings (cf. [AB], [AL], [BF], [BM], [BW], [GTZ], [KC], 
[K-RK], [M\"ol], [NS1], [NS2], [Pan], [Pau], [Zac]). To go into a 
little more detail, we mainly refer to [AL] for a general theory of 
Gr\"obner bases over rings.  Let $A=R[x_1,\ldots x_n]$ be the 
polynomial ring in $n$ variables over a commutative Noetherian ring 
$R$ with the multiplicative identity 1. Then $A$ is a Noetherian 
ring and $A$ is also a free $R$-module with the standard free 
$R$-basis $\B =\{ x^{\alpha}=x_1^{\alpha_1}\cdots 
x_n^{\alpha_n}~|~\alpha =(\alpha_1,\ldots 
,\alpha_n)\in\mathbb{N}^n\}$, where $\mathbb{N}$ denotes the set of 
all nonnegative integers. It follows from ([AL], Chapter 4) that if 
linear equations are solvable in $R$, then a Gr\"obner basis theory 
holds true for $A$ and  (in principle) finite Gr\"obner bases in the 
sense of ([AL], Definition 4.1.13) may be constructed by means of an 
analogue of Buchberger's algorithm ([AL], Algorithm 4.2.1). However, 
as one may see from the literature, except for some specific rings 
such as $R=\mathbb{Z}$ which is the ring of integers,  and $R=K[y]$ 
which is the polynomial ring in one variable $y$ over a field $K$, 
the practical implementation of such an algorithm seems to be  much 
restricted  by the coefficient ring $R$, in particular, it is rather 
annoying when $R$ has divisors of zero, or when modulo operation has 
to be considered in dealing with the operations of elements in $R$. 
To remedy such problems, in this paper we propose a ``pull-back 
method" whenever $R$ is a quotient ring of a ring $\Lambda$ over 
which Gr\"obner bases may be effectively (or even more effectively) 
constructed by means of well-implemented algorithms. More precisely, 
in Section 2 we  show that if $\Lambda$ is an {\it arbitrary} 
commutative Noetherian ring in which linear equations are solvable, 
then linear equations are solvable in any quotient ring $R=\Lambda 
/I$ of $\Lambda$, all basic results concerning  Gr\"obner bases as 
presented in ([AL], Sections 4.1 -- 4.3) hold true over  $R$ and can 
be realized over $\Lambda$. As applications of Section 2, in Section 
3 we show that if $\Lambda$ is a PID and $R=\Lambda /I$ is any 
quotient ring of $\Lambda$, then the theory of strong Gr\"obner 
bases holds true over $R$ and can be realized over $\Lambda$; while 
section 4 deals with the case where $R=E[a_1,\ldots ,a_m]$ is a 
finitely generated $E$-algebra over a commutative Noetherian ring 
$E$, and we show, by employing a nice result of [AB], that if linear 
equations are solvable in $E$ then the Gr\"obner basis theory ([AL], 
Sections 4.1 -- 4.3) holds true over $R$ and can be realized over 
$E$. In Section 5 we illustrate the practical effectiveness of 
Theorem 4.4 by focusing on a finitely generated $D$-algebra $R 
=D[\vartheta_1,\ldots ,\vartheta_m]$ where $D$ is a PID including 
the case of $D$ being a field. Finally in Section 6 we give an 
comprehensive application of the previous sections to the case where 
$R$ is a Galois ring.  \v5

Let $\Lambda$ be a commutative Noetherian ring. We now recall from
[AL] some basics on  Gr\"obner bases for ideals in the polynomial
ring $T=\Lambda [x_1,\ldots x_n]$. Let $\B =\{
x^{\alpha}=x_1^{\alpha_1}\cdots x_n^{\alpha_n}~|~\alpha
=(\alpha_1,\ldots ,\alpha_n)\in\mathbb{N}^n\}$ be the standard free
$\Lambda$-basis of $T$. For convenience, we write
$\Lambda^*=\Lambda-\{ 0\}$,  and, as usual we call $x^{\alpha}$ a
{\it monomial} of $T$. Given a monomial ordering $\prec$ on $\B$,
any nonzero element $f=\sum^m_{i=1}\lambda_ix^{\alpha (i)}$ with
$\lambda_i\in \Lambda^*$, $\alpha (i)=(\alpha_{i_1},\ldots
,\alpha_{i_n})$ and $x^{\alpha_(m)}\prec x^{\alpha
(m-1)}\prec\cdots\prec x^{\alpha (1)}$, has the associated
$$\begin{array}{l} \hbox{leading monomial:} ~\LM (f)=x^{\alpha (1)},\\
\hbox{leading coefficient:} ~\LC (f)=\lambda_1,~ and\\
\hbox{leading term:} ~\LT (f)=\LC (f)\LM (f)=\lambda_1x^{\alpha
(1)}.\end{array}$$ If $S$ is a subset of $T$, we write $\langle
S\rangle$ for the ideal of $T$ generated by $S$, and we write $\LT
(S)=\{ \LT (f)~|~f\in S\}$ for the set of leading terms of
$S$.\par

In order to have an effective division algorithm and to make
Gr\"obner bases computable over $\Lambda$, it is necessary to have
the following {\parindent=0pt\v5

{\bf 1.1. Definition} We will say that linear equations are solvable
in $\Lambda$ provided that\par

(i) Given $\lambda ,\lambda_1,\ldots ,\lambda_m\in \Lambda$, there
is an algorithm to determine whether
$\lambda\in\langle\lambda_1,\ldots ,\lambda_m\rangle$ and if it is,
to compute $\mu_1,\ldots ,\mu_m\in \Lambda$ such that $\lambda
=\lambda_1\mu_1+\cdots +\lambda_m\mu_m$;\par

(ii) Given $\lambda_1,\ldots ,\lambda_m\in \Lambda$, there is an
algorithm that computes a set of generators for the $\Lambda$-module
$$\hbox{Syz}_{\Lambda}(\lambda_1,\ldots ,\lambda_m)=\{ (\mu_1,\ldots ,\mu_m)\in \Lambda^m~|~\lambda_1\mu_1
+\cdots +\lambda_m\mu_m=0\}.$$} \par

Given $f\in T$ and a finite set $F=\{ f_1,\ldots ,f_s\}\subset T-\{
0\}$, if, under the assumption of Definition 1.1(i), there are
$\lambda_1,\ldots ,\lambda_{s'}\in \Lambda^*$, $x^{\alpha
(1)},\ldots ,x^{\alpha (s')}\in\B$, and $f_{i_1},\ldots
,f_{i_{s'}}\in F$, $1\le s'\le s$, such that
$$\begin{array}{l} \LM (f)=x^{\alpha (j)}\LM (f_{i_j}),~1\le j\le s',\\
\LT (f)=\lambda_1x^{\alpha (1)}\LT (f_{i_1})+\cdots
+\lambda_{s'}x^{\alpha (s')}\LT (f_{i_{s'}}),\end{array}$$ then $f$
can be expressed as
$$f=(\lambda_1x^{\alpha (1)}f_{i_1}+\cdots +\lambda_{s'}x^{\alpha (s')}f_{i_{s'}})+h,~h\in T~\hbox{with}~
\LM (h)\prec\LM (f).$$ In this case, $f$ is said to be {\it reduced
to} $h$ modulo $F$ in one step and is denoted by
$f~\mapright{}{}~h$. If $f\ne 0$ and $f$ cannot be reduced as above,
then we say that $f$ is {\it minimal  with respect to} $F$. Thus,
under the assumption that linear equations are solvable in $E$,
there is an effective division algorithm ([AL], Algorithm 4.1.1)
that produces elements $h_1,\ldots ,h_s,r\in T$ such that
$$f=h_1f_1+\cdots +h_sf_s+r,$$
where $\LM (h_i)\LM (f_i)\preceq\LM (f)$ for all $h_i\ne 0$, and
either $r=0$ or $r$ is minimal with $\LM (r)\preceq\LM (f)$.\par

By using the division algorithm described above, the following
theorem is proved in [AL]. But note that we modified below the
condition (ii) of ([AL], Theorem 4.1.12) in an equivalent statement.
{\parindent=0pt\v5

{\bf 1.2. Theorem}  Let $I$ be an ideal of $T$ and $\G =\{
g_1,\ldots ,g_s\}$ a finite subset of $I-\{ 0\}$. The following
statements are equivalent.\par

(i) $\langle \LT (I)\rangle =\langle \LT (\G )\rangle$.\par

(ii) If $0\ne f\in I$, then there are $\lambda_1,\ldots
,\lambda_k\in \Lambda^*$, $x^{\alpha (1)},\ldots ,x^{\alpha
(k)}\in\B$, and $g_{i_1},\ldots ,g_{i_k}\in\G$, such that
$$\begin{array}{l} \LM (f)=x^{\alpha (j)}\LM (g_{i_j}),~1\le j\le k,\\
\LT (f)=\lambda_1x^{\alpha (1)}\LT (g_{i_1})+\lambda_2x^{\alpha
(2)}\LT (g_{i_2})+\cdots + \lambda_kx^{\alpha (k)}\LT
(g_{i_k}).\end{array}$$
\par

(iii) If $0\ne f\in I$, then $f$ has a representation
$$\begin{array}{rcl} f&=&\sum^s_{i=1}h_ig_i,~h_i\in T,~g_i\in\G,\\
&{~}&\hbox{satisfying}~\LM (h_i)\LM (g_i)\preceq\LM (f)~\hbox{for
all}~h_i\ne 0.\end{array}$$\v5

{\bf 1.3. Definition} Let $I$ be an ideal of $T$ and $\G =\{
g_1,\ldots ,g_s\}$ a finite subset of $I-\{ 0\}$. If $\G$ satisfies
one of the three equivalent conditions of Theorem 1.2, then we call
$\G$ a {\it Gr\"obner basis} of $I$.}\v5

Clearly, if $\G$ is a Gr\"obner basis of the ideal $I$ then it is a
generating set for $I$. So, from now on in the text if we say that a
finite subset $\G$ of nonzero elements is a Gr\"obner basis in $T$,
it is meant that $\G$ is a Gr\"obner basis for the ideal $I=\langle
\G\rangle$. Moreover, as usual a represention of $f$ given in
Theorem 1.2(iii) is called a {\it Gr\"obner representation} of $f$
by $\G$.\par

It follows from ([AL], Chapter 4) that if linear equations are
solvable in $\Lambda$, then Gr\"obner bases in $T$ are computable,
that is, with a given subset $S=\{ f_1,\ldots ,f_t\}\subset T$,
there is an analogue of Buchberger's algorithm ([AL], Algorithm
4.2.1, Algorithm 4.2.2) for computing a finite Gr\"obner basis for
the ideal $I=\langle S\rangle$. In this case, we say that {\it
Gr\"obner bases in T are computable in the sense of} [AL]. \v5

\section*{2. The General Case: $R=\Lambda /I$}
Let $T=\Lambda [x_1,\ldots x_n]$ be the polynomial ring in $n$
variables over a commutative Noetherian ring $\Lambda$, and let $I$
be a proper  ideal of $\Lambda$, $R=\Lambda /I$. Suppose that
Gr\"obner bases in $T$ are computable in the sense of [AL]. In this
section, we show that Gr\"obner bases for ideals of the polynomial
ring $A=R[x_1,\ldots ,x_n]$ are well defined and can be obtained by
constructing Gr\"obner bases in $T$. Notations and conventions fixed
in the last section are maintained. \v5

Since Gr\"obner bases in $T$ are computable by the assumption,
linear equations are solvable in $\Lambda$. We first show that
linear equations are solvable in $R=\Lambda/I$, thereby (in
principle) Gr\"obner bases in $A$ are computable. As usual, if
$\lambda\in \Lambda $ then we write $\OV{\lambda}$ for the coset
$\lambda +I$ in $R$ represented by $\lambda$.{\parindent=0pt\v5

{\bf 2.1. Proposition} Suppose that $I=\langle \nu_1,\ldots
,\nu_s\rangle$ where $\nu_1,\ldots ,\nu_s\in I$, and that linear
equations are solvable in $\Lambda$. Then linear equations are
solvable in $R=\Lambda /I$.\vskip 6pt

{\bf Proof} We show that $R$ satisfies the two conditions of
Definition 1.1.}\par

(i) Given $\OV{\lambda} ,\OV{\lambda}_1,\ldots ,\OV{\lambda}_m\in R$
represented by $\lambda ,\lambda_1,\ldots ,\lambda_m\in \Lambda$
respectively, since $I=\langle \nu_1,\ldots ,\nu_s\rangle$ we have
$$\begin{array}{rcl} \OV{\lambda}\in\langle\OV{\lambda}_1,\ldots ,\OV{\lambda}_m\rangle&\Leftrightarrow&
\hbox{there are}~\OV{\mu}_1,\ldots ,\OV{\mu}_m\in R~\hbox{such
that}~\OV{\lambda}=\sum^m_{i=1}\OV{\lambda}_i\OV{\mu}_i\\
&\Leftrightarrow&\lambda-\sum^m_{i=1}\lambda_i\mu_i\in I=\langle
\nu_1,\ldots ,\nu_s\rangle\\
&\Leftrightarrow&\lambda\in J=\langle\lambda_1,\ldots
,\lambda_m,\nu_1,\ldots ,\nu_s\rangle\subset \Lambda.\end{array}$$
By the assumption, there is an algorithm to determine whether
$\lambda\in J$ and if it is, to compute $\mu_1,\ldots
,\mu_m,\xi_1,\ldots ,\xi_s\in \Lambda$ such that $\lambda
=\sum^m_{i=1}\lambda_i\mu_i+\sum^s_{j=1}\nu_j\xi_j$, i.e.,
$\OV{\lambda}=\sum^m_{i=1}\OV{\lambda}_i\OV{\mu}_i$, or
equivalently, $\OV{\lambda}\in\langle\OV{\lambda}_1,\ldots
,\OV{\lambda}_m\rangle$.\par

(ii) Given $\OV{\lambda}_1,\ldots ,\OV{\lambda}_m\in R$ represented
by $\lambda_1,\ldots ,\lambda_m\in \Lambda$ respectively, since
$I=\langle \nu_1,\ldots ,\nu_s\rangle$ we have
$$\begin{array}{l} (\OV{\mu}_1,\ldots ,\OV{\mu}_m)\in~\hbox{Syz}_R(\OV{\lambda}_1,\ldots ,\OV{\lambda}_m)\subset R^m\\
\quad\quad\quad\quad\quad\quad\quad \Leftrightarrow\sum^m_{i=1}\OV{\lambda}_i\OV{\mu}_i=0\\
\quad\quad\quad\quad\quad\quad\quad \Leftrightarrow\sum^m_{i=1}\lambda_i\mu_i\in I\\
\quad\quad\quad\quad\quad\quad\quad \Leftrightarrow\hbox{there
are}~\xi_1,\ldots ,\xi_s\in \Lambda~\hbox{such
that}~\sum^m_{i=1}\lambda_i\mu_i+\sum^s_{j=1}\nu_j\xi_j=0\\
\quad\quad\quad\quad\quad\quad\quad \Leftrightarrow(\mu_1,\ldots
,\mu_m,\xi_1,\ldots ,\xi_s)\in~\hbox{Syz}_E(\lambda_1,\ldots
,\lambda_m,\nu_1,\ldots ,\nu_s)\subset \Lambda^{m+s}.\end{array}$$
By the assumption, there is an algorithm that computes a set of
generators $\{ V_1,\ldots V_q\}$ for the $\Lambda$-module
$\hbox{Syz}_{\Lambda}(\lambda_1,\ldots ,\lambda_m,\nu_1,\ldots
,\nu_s)\subset E^{m+s}$. Consider the $\Lambda$-module epimorphism
$$\begin{array}{cccc} \psi :&\Lambda^{m+s}&\mapright{}{}&R^m\\
&(\mu_1,\ldots ,\mu_m,\xi_1,\ldots
,\xi_s)&\mapsto&(\OV{\mu}_1,\ldots ,\OV{\mu}_m)\end{array}$$ It
follows from the above discussion that
$$\psi\left (\hbox{Syz}_{\Lambda}(\lambda_1,\ldots ,\lambda_m,\nu_1,\ldots ,\nu_s)\right )=~\hbox{Syz}_R(\OV{\lambda}_1,
\ldots ,\OV{\lambda}_m).$$ Since the $\Lambda$-action of $\Lambda$
on $R^m$ and the $R$-action of $R$ on $R^m$ coincide,  a generating
set of the $\Lambda$-module Syz$_{\Lambda}(\lambda_1,\ldots
,\lambda_m,\nu_1,\ldots ,\nu_s)$ is mapped to a generating set of
the $R$-module Syz$_R(\OV{\lambda}_1,\ldots ,\OV{\lambda}_m)$,
thereby $\{ \psi (V_1),\ldots \psi (V_q)\}$ generates the $R$-module
Syz$_R(\OV{\lambda}_1,\ldots ,\OV{\lambda}_m)$.\QED \v5

Although Proposition 2.1 tells us that in principle Gr\"obner bases
in the polynomial ring $A=R[x_1,\ldots ,x_n]$ in any $n\ge 1$
variables are computable, its proof indeed inspires us to obtain
Gr\"obner bases for ideals in $A$ by constructing Gr\"obner bases in
the polynomial ring $T=\Lambda [x_1,\ldots ,x_n]$. To go further,
without confusion we let $\B =\{ x^{\alpha}=x_1^{\alpha_1}\cdots
x_n^{\alpha_n}~|~\alpha =(\alpha_1,\ldots
,\alpha_n)\in\mathbb{N}^n\}$ denote the standard free
$\Lambda$-basis of $T$ and the standard free $R$-basis of $A$, and
we fix a monomial ordering $\prec$ on $\B$.  Let $\varphi$: $T\r A$
be the canonical ring epimorphism with $\varphi
(\sum_i\lambda_ix^{\alpha (i)})=\sum_i\OV{\lambda}_ix^{\alpha (i)}$,
where $\lambda_i\in \Lambda$, $x^{\alpha (i)}\in\B$, and
$\OV{\lambda}_i=\lambda_i+I$ is the coset in $R=\Lambda /I$
represented by $\lambda_i$. Recalling the definition of the leading
term given in Section 1, the lemma given below is
clear.{\parindent=0pt\v5

{\bf 2.2. Lemma} (i) If $0\ne f=\sum^m_{i=1}\lambda_ix^{\alpha
(i)}\in T$ which has the leading term $\LT (f)=\lambda_1x^{\alpha
(1)}$, then $\varphi (f)\ne 0$ and $\LT (\varphi (f))=\varphi (\LT
(f))$ if and only if $\lambda_1\not\in I$.\par

(ii) If $0\ne \OV f=\sum_{i=1}^m\OV{\lambda}_ix^{\alpha (i)}\in A$
which has $\LT (\OV f)=\OV{\lambda}_1x^{\alpha (1)}$, then, putting
$f=\sum^m_{i=1}\lambda_ix^{\alpha (i)}\in T$, we have $\LT
(f)=\lambda_1x^{\alpha (1)}$ and $\LT (\varphi (f))=\varphi (\LT
(f))=\LT (\OV f)$.\par \QED\v5

{\bf 2.3. Theorem} Let the rings $T=\Lambda [x_1,\ldots ,x_n]$,
$R=\Lambda /I$, $A=R[x_1,\ldots ,x_n]$, the canonical ring
epimorphism $\varphi$: $T\r A$, and all notations be as above.
Suppose that the ideal $I$ is generated by $\nu_1,\ldots ,\nu_s\in
\Lambda$, and let $\OV J=\langle\OV S\rangle$ be an ideal of $A$
generated by the set $\OV S=\{\OV f_1,\ldots ,\OV f_t\}$, where $\OV
f_i=\sum_j\OV{\lambda}_{ij}x^{\alpha (j)}$ with $\lambda_{ij}\in
\Lambda$, $1\le i\le t$. Consider in $T$ the set of elements $S=\{
f_1,\ldots ,f_t,\nu_1,\ldots ,\nu_s\}$, where
$f_i=\sum_j\lambda_{ij}x^{\alpha (j)}$, and let $J=\langle S\rangle$
be the ideal of $T$ generated by $S$. The following statements
hold.\par

(i) If, with respect to the given monomial ordering $\prec$, $\G =\{
g_1,\ldots ,g_m\}$ is a Gr\"obner basis of $J$ constructed by means
of ([AL], Algorithm 4.2.1) using the initial input-data $S$ in $T$,
then for each $g_i\not\in\{ \nu_1,\ldots ,\nu_s\}$, we have $\varphi
(g_i)\ne 0$, $\LT (\varphi (g_i))=\varphi (\LT (g_i))$ and hence
$\LM (\varphi (g_i))=\LM (g_i)$.\par

(ii) Let $\G$ be the Gr\"obner basis of $J$ presented in (i).
Considering the image $\varphi (\G )$ in $A$, the set
$$\OV{\G}=\{ \OV g_k=\varphi (g_k)~|~g_k\in\G ,~g_k\ne \nu_j,~1\le j\le s\}=\varphi (\G )-\{ 0\}$$
is a Gr\"obner basis for $\OV J$ with respect to the monomial
ordering $\prec$.\par

(iii) Let $\G$ be the Gr\"obner basis of $J$ presented in (i). Then
each $g_k\in\G-\{\nu_1,\ldots ,\nu_s\}$ has a representation
$g_k=g_k'+g_k''$ in which
$$\begin{array}{l} g_k'=\sum_q\mu_{kq}x^{\alpha (q)}~\hbox{with}~\mu_{kq}\not\in I~\hbox{for all}~q,\\
g_k''=\sum_{\ell}\mu_{k\ell}x^{\alpha (\ell
)}~\hbox{with}~\mu_{k\ell}\in I~\hbox{for all}~\ell,\end{array}$$
and such a representation can be algorithmically determined. Hence
each element of the Gr\"obner basis $\OV{\G}$ obtained in (ii) has a
``real representation" in $A$, i.e., $\OV g_k=\varphi (g_k)=\varphi
(g_k')=\sum_q\OV{\mu}_{kq}x^{\alpha (q)}$ with all the
$\OV{\mu}_{kq}\ne 0$. \vskip 6pt

{\bf Proof} (i) Let $g_i\in\G$ with $g_i\not\in\{ \nu_1,\ldots
,\nu_q\}$. If $g_i\in\{ f_1,\ldots ,f_t\}$, then $g_i$ has the
desired property by Lemma 2.2. If $g_i\not\in\{ f_1,\ldots
,f_t,\nu_1,\ldots ,\nu_s\}$, then it follows from ([AL], Algorithm
4.2.1) that $g_i$ is obtained by passing through some $i$-th round
executing the While loop, that is, $g_i:=r$, where $r$ appears as
the remainder of a reduction
$$h_1g_1+\cdots +h_{\ell}g_{\ell}~\mapright{G'}{}_+r$$
which is minimal with respect to $G'$. Note that $S\subset G'$ since
$S$ is the initial input-data running the algorithm. So, $g_i$
cannot be further reduced modulo $G'$, in particular, it cannot be
reduced modulo $\{ \nu_1,\ldots ,\nu_s\}$, i.e., if $\LT
(g_i)=\lambda x^{\alpha}$, then there do not exist $\lambda_1,\ldots
,\lambda_s\in \Lambda$ such that
$$\lambda x^{\alpha}=(\lambda_1\nu_1+\cdots +\lambda_s\nu_s)x^{\alpha}.$$
This shows that if we write $g_i=\lambda
x^{\alpha}+\sum_j\lambda_jx^{\alpha (j)}$ with $x^{\alpha (j)}\prec
x^{\alpha}$, then $\varphi
(g_i)=\OV{\lambda}x^{\alpha}+\sum_j\OV{\lambda}_jx^{\alpha (j)}$
with $\OV{\lambda}\ne 0$ and hence $\LT (\varphi
(g_i))=\OV{\lambda}x^{\alpha}=\varphi (\LT (g_i))$, $\LM (\varphi
(g_i))=x^{\alpha}=\LM (g_i)$.}\par

(ii) Since $J=\varphi^{-1}(\OV J)=\{ f\in T~|~\varphi (f)\in \OV
J\}$ and hence $\varphi (J)=\OV J$, we first note that
$\OV{\G}\subset \OV J$. For any $0\ne \OV
f=\sum^m_{i=1}\OV{\lambda}_ix^{\alpha (i)}\in\OV J$ with $\LT (\OV
f)=\OV{\lambda}_1x^{\alpha (1)}$, putting
$f=\sum^m_{i=1}\lambda_ix^{\alpha (i)}$ in $T$ we have $f\in J$ with
$\LT (f)=\lambda_1x^{\alpha (1)}$ and, by Lemma 2.2,
$$\LT (\varphi (f))=\varphi (\LT (f))=\LT (\OV f)=\OV{\lambda}_1x^{\alpha (1)}. \eqno{(1)}$$
Moreover, by Theorem 1.2, there are $\mu_1,\ldots ,\mu_k\in
\Lambda^*=\Lambda -\{ 0\}$, $x^{\beta (1)},\ldots ,x^{\beta
(k)}\in\B$, where $\beta (j)=(\beta_{j1},\ldots
,\beta_{jn})\in\mathbb{N}^n$, $1\le j\le k$, and $g_{i_1},\ldots
,g_{i_k}\in\G$ such that
$$\begin{array}{l} x^{\alpha (1)}=\LM (f)=x^{\beta (j)}\LM (g_{i_j}),~1\le j\le k,\\
\lambda_1x^{\alpha (1)}=\LT (f)=\mu_1x^{\beta (1)}\LT
(g_{i_1})+\cdots +\mu_kx^{\beta (k)}\LT
(g_{i_k}).\end{array}\eqno{(2)}$$ Since $\LT (\OV
f)=\OV{\lambda}_1x^{\alpha (1)}\ne 0$, by the definition of
$\varphi$ we may assume that $\mu_j\not\in I$ and $g_{i_j}\not\in\{
\nu_1,\ldots ,\nu_s\}$, $1\le j\le k$. Combining (1) and (2), it
follows from (i) that
$$\begin{array}{l}\hskip .2truecm\LM (\OV f)=x^{\alpha (1)}=\LM (f)=x^{\beta (j)}\LM (\varphi (g_{i_j})),~1\le j\le k,\\
\begin{array}{rcl}
\LT (\OV f)=\OV{\lambda}_1x^{\alpha (1)}=\varphi (\LT (f))&=&
\OV{\mu}_1x^{\beta (1)}\varphi (\LT (g_{i_1}))+\cdots
+\OV{\mu}_kx^{\beta (k)}\varphi (\LT (g_{i_k}))\\
&=&\OV{\mu}_1x^{\beta (1)}\LT (\varphi (g_{i_1}))+\cdots
+\OV{\mu}_kx^{\beta (1)}\LT (\varphi
(g_{i_k})).\end{array}\end{array}$$ This shows that $\OV{\G}$
satisfies the condition Theorem 1.2(ii), thereby $\OV{\G}$ is a
Gr\"obner basis for $\OV J$.\par

(iii) Since $\G$ is finite and since linear equations are solvable
in $\Lambda$ by our assumption, it follows that there is an
algorithm to determine whether a coefficient of $g_k$ is in
$I=\langle\nu_1,\ldots ,\nu_s\rangle$.\par \QED\v5

In ([AL], Section 4.3) several basic applications of Gr\"obner bases
over rings are presented. We next show that those basic applications
of Grobner bases with coefficients in the quotient ring $R=\Lambda
/I$ may be realized by using Gr\"obner bases with coefficients in
the ring $\Lambda$. In what follows, we let the rings $T=\Lambda
[x_1,\ldots ,x_n]$, $R=\Lambda /I$ with $I=\langle \nu_1,\ldots
\nu_s\rangle$, $A=R[x_1,\ldots ,x_n]$, the canonical ring 
epimorphism $\varphi$: $T\r A$, the ideal $\OV J=\langle\OV 
S\rangle$ of $A$, and the ideal $J=\langle S\rangle$ of $T$ be as in 
Theorem 2.3. So we have the Gr\"obner basis $\G =\{ g_1,\ldots 
,g_m\}$ for $J$ and the Gr\"obner basis $\OV{\G}=\varphi (\G )-\{ 
0\}$ for $\OV J$. {\parindent=0pt\v5

{\bf 2.4. Proposition} (Ideal membership problem) Let $0\ne \OV
f=\sum_{i=1}^s\OV{\lambda}_ix^{\alpha (i)}\in A$ with all the
$\OV{\lambda}_i\ne 0$. Then $\OV f\in \OV J$ if and only if
$f=\sum_{i=1}^s\lambda_ix^{\alpha (i)}\in J$ if and only if
$f~\mapright{\G}{}_+~0$, i.e., by the division by $\G$ ([AL],
Algorithm 4.1.1), $f$ has a Gr\"obner representation:
$$\begin{array}{rcl} f&=&\sum^m_{i=1}h_ig_i,~h_i\in T,~g_i\in\G,\\
&{~}&\hbox{satisfying}~\LM (h_i)\LM (g_i)\preceq\LM (f)~\hbox{for
all}~h_i\ne 0.\end{array}$$\par

{\bf Proof} Since $J=\varphi^{-1}(\OV J)$ and $\varphi (f)=\OV f$,
the assertion is clear by Theorem 2.3.\QED\v5

{\bf Remark} It follows from Theorem 2.3 and Proposition 2.4 that if
$0\ne \OV f\in \OV J$, then the Gr\"obner representation of $f$ by
$\G$ gives rise to a Gr\"obner representation of $\OV f$ by
$\OV{\G}$:
$$\begin{array}{rcl} \OV f=\varphi (f)&=&\sum^m_{i=1}\varphi (h_i)\varphi (g_i),~h_i\in T,~g_i\in\G,\\
&=&\sum_k\varphi (h_k)\OV g_k,~0\ne \varphi (g_k)=\OV g_k\in\OV{\G}\\
&{~}&\hbox{satisfying}~\LM (\varphi (h_i))\LM (\OV g_k)\preceq\LM
(\OV f)~\hbox{for all}~\varphi (h_k)\ne 0.\end{array}$$\par

{\bf 2.5. Proposition} A complete set of coset representatives for
$A/\OV J$ can be obtained by constructing a complete set of coset
representatives for $T/J$ as in ([AL], Theorem 4.3.3) provided
linear equations are solvable in $\Lambda$ and $\Lambda$ has
effective coset representatives (in the sense of ([AL],
P.226).\vskip 6pt

{\bf Proof} This assertion  follows from the fact that
$J=\varphi^{-1}(\OV J)$ and hence $A/\OV J\cong T/J$. \QED\v5

{\bf 2.6. Theorem} With $T=\Lambda [x_1,\ldots ,x_n,y_1,\ldots
,y_m]$, $A=R[x_1,\ldots ,x_n,y_1,\ldots ,y_m]$, where $R=\Lambda /I$
is as before, if $\G$ is a Gr\"obner basis for $J$ with respect to
an elimination ordering $\prec$ on the standard free $\Lambda$-basis
$\B$ of $T$ with the $y$ variables larger than the $x$ variables,
then $\OV{\G}\cap R[x_1,\ldots ,x_n]$ is a Gr\"obner basis for $\OV
J\cap R[x_1,\ldots ,x_n]$ with respect to the $\prec$ restricted on
$R[x_1,\ldots ,x_n]$.\vskip 6pt

{\bf Proof} By Theorem 2.3,  $\OV{\G}=\varphi (\G )-\{ 0\}$ is a
Gr\"obner basis for $\OV J$ with respect to the same elimination
ordering on the standard free $R$-basis $\B$ of $A$, this assertion
follows from ([AL], Theorem 4.3.6).\QED\v5

{\bf 2.7. Corollary} With $I=\langle \nu_1,\ldots ,\nu_s\rangle$ and
$R=\Lambda /I$, it follows from  Theorem 2.6 and ([AL], Proposition
4.3.9, Proposition 4.3.11, Theorem 4.3.13) that the following
results hold.\par

(i) A generating set (indeed a Gr\"obner basis) for $\OV J_1\cap\OV
J_2$, where $\OV J_1=\langle \OV f_1,\ldots ,\OV f_{m_1}\rangle$ and
$\OV J_2=\langle \OV h_1,\ldots ,\OV h_{m_2}\rangle$ are ideals of
$A=R[x_1,\ldots ,x_n]$, can be obtained by constructing a Gr\"obner
basis $\G$ for the ideal $J=\langle yf_1,\ldots
,yf_{m_1},(1-y)h_1,\ldots ,(1-y)h_{m_2},\nu_1,\ldots ,\nu_s\rangle$
in the polynomial $\Lambda$-algebra $\Lambda [y, x_1,\ldots ,x_n]$.
\par

(ii) Let $\OV J_1$ and $\OV J_2$ be as in (i). Then the ideal
quotient
$$\OV J_1:\OV J_2=\{ \OV f\in A~|~\OV f\OV J_2\subseteq\OV J_1\}=\cap^{m_2}_{j=1}\OV J_1:\langle \OV h_j\rangle$$
can be obtained by constructing a generating set (indeed a Gr\"obner
basis) of $\OV J_1\cap\langle \OV h_j\rangle$ as in (i).\par

(iii) If $\psi$: $B\r A$ is an $R$-algebra homomorphism of the
polynomial $R$-algebra $B=R[y_1,\ldots ,y_m]$ to the polynomial
$R$-algebra $A=R[x_1,\ldots ,x_n]$ such that $\psi (y_i)=\OV f_i$,
$1\le i\le m$ then the kernel Ker$\psi$ of $\psi$ can be obtained by
constructing a Gr\"obner basis $\G$ for the ideal $J=\langle
y_1-f_1,\ldots ,y_m-f_m,\nu_1,\ldots ,\nu_s\rangle$ in the
polynomial $\Lambda$-algebra $\Lambda [y_1,\ldots ,y_m,x_1,\ldots
,x_n]$.\par\QED}\v5

Finally, with $I=\langle \nu_1,\ldots ,\nu_s\rangle$ and $R=\Lambda
/I$, we demonstrate how to obtain a generating set for the syzygy
module Syz$_A(\OV f_1,\ldots ,\OV f_t)$ of a set of nonzero elements
$\{ \OV f_1,\ldots ,\OV f_t\}\subset A=R[x_1,\ldots ,x_n]$ by
working out  necessary data in  $T=\Lambda [x_1,\ldots
,x_n]$.{\parindent=0pt\v5

{\bf 2.8. Theorem} Let $\varphi$: $T\r A$ be the canonical ring
epimorphism as before. A generating set for the syzygy module
Syz$_A(\OV f_1,\ldots ,\OV f_t)$ of a set of nonzero elements $\OV
S=\{ \OV f_1,\ldots ,\OV f_t\}\subset A=R[x_1,\ldots ,x_n]$ can be
obtained as follows.\par

{\bf Step 1.} Starting with the initial input-data $S=\{f_1,\ldots
,f_t,\nu_1,\ldots ,\nu_s\}$, ([AL], Algorithm 4.2.2), which is the
Gr\"obner basis algorithm using M\"oller's technique [M\"ol],
produces a Gr\"obner basis $\G =\{ g_1,\ldots g_m\}$ for the ideal
$J=\langle S\rangle$ of $T$ and, at the same time, outputs a
homogeneous generating set $H=\{ V_1,\ldots ,V_{\ell}\}$ of the
syzygy module Syz$_T(\LT (g_1),\ldots ,\LT (g_m))\subset T^m$.\par

{\bf Step 2.} By Theorem 2.3, $\OV{\G}=\varphi (\G )-\{ 0\}$ is a
Gr\"obner basis for the ideal $\OV J=\langle\OV S\rangle$ of $A$.
Suppose that $\OV{\G }=\{\OV g_{i_1},\ldots , \OV g_{i_d}\}$. We
conclude that the $H$ obtained in Step 1 gives rise to a homogeneous
generating set $\OV H$ for the syzygy module Syz$_T(\LT (\OV
g_{i_1}),\ldots ,\LT (\OV g_{i_d}))$. More precisely,  rewriting $\G
=\{ g_{i_1},\ldots ,g_{i_d},\nu_1,\ldots ,\nu_s\}$, we have $\OV
H=\phi (H)$, where $\phi$ is the $T$-module epimorphism
$$\begin{array}{cccc} \phi :&T^m&\mapright{}{}&A^d\\
&(h_1,\ldots ,h_d,h_1',\ldots ,h_s')&\mapsto&(\OV h_1,\ldots ,\OV
h_d)\end{array}$$\par

{\bf Step 3.} With $\OV{\G}$ and $\OV H$ obtained in Step 2,  a
generating set for the syzygy module Syz$_A(\OV f_1,\ldots ,\OV
f_t)$ can then be obtained by ([AL], Theorem 4.3.16).\vskip 6pt

{\bf Proof} We need only to prove the conclusion on $\OV H$ in Step
2. To this end, rewrite $\G =\{ g_{i_1},\ldots ,g_{i_d},\nu_1,\ldots
,\nu_s\}$. For $\OV h_1,\ldots ,\OV h_d\in A$ represented by
$h_1,\ldots ,h_d\in T$ respectively, since Ker$\varphi$ is the ideal
generated by $\nu_1,\ldots ,\nu_q$ in $T$, we have by Theorem 2.3
that
$$\begin{array}{l} (\OV h_1,\ldots ,\OV h_d)\in~\hbox{Syz}_A(\LT (\OV g_{i_1}),\ldots ,\LT(\OV g_{i_d}))\subset A^d\\
\quad\quad\quad\quad\quad\quad\quad \Leftrightarrow \sum^d_{j=1}\OV
h_j\LT (\OV g_{i_j})=0\\
\quad\quad\quad\quad\quad\quad\quad \Leftrightarrow \varphi
(\sum^d_{j=1}h_j\LT (g_{i_j}))=\sum^d_{j=1}\varphi (h_j)\varphi (\LT
(g_{i_j})=\sum^d_{j=1}\varphi (h_j)\LT (\varphi (g_{i_j}))=0\\
\quad\quad\quad\quad\quad\quad\quad \Leftrightarrow
\sum^d_{j=1}h_j\LT (g_{i_j})\in~\hbox{Ker}\varphi\\
\quad\quad\quad\quad\quad\quad\quad \Leftrightarrow \hbox{there
are}~h_1',\ldots ,h_s'\in T~\hbox{such that}~\sum^d_{j=1}h_j\LT
(g_{i_j})+\sum^s_{k=1}h_k'\nu_k=0\\
\quad\quad\quad\quad\quad\quad\quad \Leftrightarrow (h_1,\ldots
,h_d,h_1',\ldots ,h_s')\in~\hbox{Syz}_T(\LT (g_1),\ldots ,\LT
(g_m))\in T^m\end{array}$$
Consider the $T$-module epimorphism $\phi$ presented above. It
follows from the above discussion that
$$\phi\left (\hbox{Syz}_T(\LT (g_{i_1}),\ldots ,\LT (g_{i_d}),\nu_1,\ldots ,\nu_s)\right )
=~\hbox{Syz}_A(\LT (\OV g_{i_1}), \ldots ,\LT (\OV g_{i_d})).$$
Since the $T$-action of $T$ on $A^m$ and the $A$-action of $A$ on
$A^d$ coincide,  the homogeneous generating set $H$ of the
$T$-module Syz$_T(\LT (g_{i_1}),\ldots ,\LT (g_{i_d}),\nu_1,\ldots
,\nu_s)$ is mapped to a homogeneous generating set $\OV H$ of the
$A$-module Syz$_A(\LT (\OV g_{i_1}), \ldots ,\LT (\OV g_{i_d})).$
}\v5

\section*{3. The Case of $R=\Lambda /I$ with $\Lambda$ a PID}
It is well known that the theory of Gr\"obner bases for a polynomial
ring $T=\Lambda [x_1,\ldots ,x_n]$ over a principal ideal domain
(PID) $\Lambda$ is much better developed, because of the fact that
division of elements in many PID's (especially in the ring of
integers $\mathbb{Z}$) can be effectively implemented on computer,
and that  a finite strong Gr\"obner basis $\G$ (see the definition
below)  can be constructed by using only $S$-polynomials for
reductions ([AL], Algorithm 4.5.1).  In this section we show that if
$I$ is a proper ideal of a PID $\Lambda$ and $R=\Lambda /I$, then
strong Gr\"obner bases for ideals of the polynomial ring
$A=R[x_1,\ldots ,x_n]$ are well defined and can be obtained by
constructing strong Gr\"obner bases in $T=\Lambda [x_1,\ldots
,x_n]$. Consequently, all basic applications of Gr\"obner bases at
the level of $A$  can be realized by using strong Gr\"obner bases at
the level of $T$, as presented in the last section. Notations and
conventions used in previous sections are maintained. \v5

We start with an {\it arbitrary} commutative Noetherian ring
$\Lambda$ in which linear equations are solvable. As before, we
write $T=\Lambda [x_1,\ldots ,x_n]$; for a given proper ideal
$I\subset\Lambda$, we write $R=\Lambda /I$,  $A=R[x_1,\ldots ,x_n]$;
and we let $\B =\{ x^{\alpha}=x_1^{\alpha_1}\cdots
x_n^{\alpha_n}~|~\alpha =(\alpha_1,\ldots
,\alpha_n)\in\mathbb{N}^n\}$ denote the standard free
$\Lambda$-basis of $T$ and the standard free $R$-basis of $A$.
Moreover, we fix a monomial ordering $\prec$ on $\B$.  Let
$\varphi$: $T\r A$ be the canonical ring epimorphism with $\varphi
(\sum_i\lambda_ix^{\alpha (i)})=\sum_i\OV{\lambda}_ix^{\alpha (i)}$,
where $\lambda_i\in \Lambda$, $x^{\alpha (i)}\in\B$, and
$\OV{\lambda}_i=\lambda_i+I$ is the coset in $R=\Lambda /I$
represented by $\lambda_i$.{\parindent=0pt\v5

{\bf 3.1. Definition} Let $J$ be an ideal of $T$ and $\G =\{
g_1,\ldots g_t\}$ a subset of nonzero elements in $J$. If, for each
$0\ne f\in J$, there exist $\lambda\in\Lambda$, $x^{\alpha}\in\B$
and  $g_i\in\G$ such that $\LT (f)=\lambda x^{\alpha}\LT (g_i)$,
then $\G$ is called a {\it strong Gr\"obner basis} for $J$.}\v5

Let $\Lambda$ be a PID  and $R=\Lambda /I$ a quotient ring of
$\Lambda$. Then it follows from Section 2 that Gr\"obner bases for
ideals of the  polynomial ring $A=R[x_1\ldots ,x_n]$, in the sense
of Definition 1.3, are well defined and can be obtained by
constructing Gr\"obner bases in $T=\Lambda [x_1,\ldots ,x_n]$.
{\parindent=0pt\v5

{\bf 3.2. Theorem} With notation as fixed above, let  the ideal
$I=\langle \nu\rangle$ be generated by $\nu$, and let $\OV
J=\langle\OV S\rangle$ be an ideal of $A$ generated by the set of
nonzero elements $\OV S=\{\OV f_1,\ldots ,\OV f_t\}$, where $\OV
f_i=\sum_j\OV{\lambda}_{ij}x^{\alpha (j)}$, $1\le i\le s$. Consider
in $T$ the set of elements $S=\{ f_1,\ldots ,f_t,\nu\}$, where
$f_i=\sum_j\lambda_{ij}x^{\alpha (j)}$, and let $J=\langle S\rangle$
be the ideal of $T$ generated by $S$.  If, with respect to a given
monomial ordering $\prec$, $\G =\{ g_1,\ldots ,g_m\}$ is a strong
Gr\"obner basis of $J$ constructed by means of ([AL], Algorithm
4.5.1) using the initial input-data $S$ in $T$, then\par

(i) the set
$$\OV{\G}=\{ \OV g_k=\varphi (g_k)~|~g_k\in\G ,~g_k\ne \nu\}=\varphi (\G )-\{ 0\}$$
is a strong Gr\"obner basis for $\OV J$ in the sense of Definition
3.1; and\par

(ii) each $g_k\in\G-\{\nu\}$ has a representation $g_k=g_k'+g_k''$
in which
$$\begin{array}{l} g_k'=\sum_q\mu_{kq}x^{\alpha (q)}~\hbox{with}~\mu_{kq}\not\in I~\hbox{for all}~q,\\
g_k''=\sum_{\ell}\mu_{k\ell}x^{\alpha (\ell
)}~\hbox{with}~\mu_{k\ell}\in I~\hbox{for all}~\ell,\end{array}$$
and such a representation can be algorithmically determined. Hence
each element of the Gr\"obner basis $\OV{\G}$ obtained in (ii) has a
``real representation" in $A$, i.e., $\OV g_k=\varphi (g_k)=\varphi
(g_k')=\sum_q\OV{\mu}_{kq}x^{\alpha (q)}$ with all the
$\OV{\mu}_{kq}\ne 0$.\vskip 6pt

{\bf Proof} (i) By ([AL], Lemma 4.5.8), the strong Gr\"obner basis
$\G$ is first of all a Gr\"obner basis for $J$ in the sense of
Definition 1.3. It follows that $\G$ satisfies Theorem 2.3(i). So,
actually as in the proof of Theorem 2.3(ii), if $0\ne \OV
f=\sum^m_{i=1}\OV{\lambda}_ix^{\alpha (i)}\in\OV J$ with $\LT (\OV
f)=\OV{\lambda}_1x^{\alpha (1)}$, then, putting
$f=\sum^m_{i=1}\lambda_ix^{\alpha (i)}$ in $T$ we have $f\in J$ with
$\LT (f)=\lambda_1x^{\alpha (1)}$. Since $\G$ is a strong Gr\"obner
basis of $J$, there exist $\lambda\in\Lambda$, $x^{\alpha}\in\B$ and
$g_i\in\G$ such that $\LT (f)=\lambda_1x^{\alpha (1)}=\lambda
x^{\alpha}\LT (g_i)$. Noticing $\OV{\lambda}_1x^{\alpha (1)}\ne 0$,
we conclude that $g_i\ne \nu$. Therefore, by Theorem 2.3(i) we have
$$\LT (\OV f)=\OV{\lambda}_1x^{\alpha (1)}=\varphi (\LT
(f))=\OV{\lambda}x^{\alpha}\varphi (\LT
(g_i))=\OV{\lambda}x^{\alpha}\LT (\varphi (g_i)).$$ This shows that
$\OV{\G}$ is a strong Gr\"obner basis for $\OV J$ in the sense of
Definition 3.1.\par

(ii) Since $\G$ is finite and since linear equations are solvable in
the PID $\Lambda$, it follows that there is an algorithm to
determine whether a coefficient of $g_k$ is in
$I=\langle\nu\rangle$.\QED}\v5

By Theorem 3.2 and Section 2, the next theorem is
straightforward.{\parindent=0pt\v5

{\bf 3.3. Theorem} Let the  Gr\"obner bases $\G$ and $\OV{\G}$ be as
in Theorem 3.2.  Then, similar results as presented in Proposition
2.4 -- Theorem 2.8 hold, that is, all basic applications of the
Gr\"obner basis $\OV{\G}$ at the level of $A=(\Lambda /I)[x_1,\ldots
,x_n]$ can be realized by using the strong Gr\"obner basis $\G$ at
the level of $T=\Lambda [x_1,\ldots ,x_n]$.\par\QED}

In view of the practical calculation, perhaps the most effective
application of Theorem 3.2 and Theorem 3.3 should be to the case
where $\Lambda =\mathbb{Z}$  and $R=\mathbb{Z}/\langle m\rangle
=\mathbb{Z}_m$ (the residue class ring modulo $m$), in particular,
when $m$ is not a prime. This advantage will be further illustrated
in Section 6.\par

Although Theorem 3.2 and Theorem 3.3 apply also to the case where
$\Lambda =K[y]$ is a polynomial ring in one variable over a field
$K$, better results will be obtained in Section 5, namely we will
see that the Gr\"obner basis $\OV{\G}$ can be obtained by
constructing a  Gr\"obner basis $\G$ in the polynomial ring
$K[y,x_1,\ldots x_n]$ by means of the classical Buchberger
algorithm, thereby all basic applications of $\OV G$ at the level of
$A$ can be realized by using $\G$ at the level of $K[y,x_1,\ldots
x_n]$. \v5

\section*{4. The Case of $R=E[a_1,\ldots ,a_m]$}
Let $E$ be a commutative Noetherian ring in which linear equations
are solvable. Then Gr\"obner bases in the polynomial ring
$E[y_1,\ldots ,y_k]$ in any $k\ge 1$ variables are computable in the
sense of [AL]. Let $R=E[a_1,\ldots ,a_m]$ be a finitely generated
$E$-algebra with generating set $\{ a_1,\ldots ,a_m\}$. Then $R\cong
E[y_1,\ldots ,y_m]/I$, where $E[y_1,\ldots ,y_m]$ is the polynomial
$E$-algebra in variables $y_1,\ldots ,y_m$, and $I$ is the kernel of
the canonical $E$-algebra epimorphism $E[y_1,\ldots ,y_m]\r R$ with
$y_i\mapsto a_i$, $1\le i\le m$.  In this section we show that
Gr\"obner bases for ideals of the polynomial ring $A=R[x_1,\ldots
,x_n]$, in the sense of Definition 1.3, are well defined, and that
if a set of generators of $I$ is known, say $I=\langle \nu_1,\ldots
\nu_s\rangle$, then Gr\"obner bases in $A$ can be obtained by
constructing Gr\"obner bases in the polynomial $E$-algebra
$T=E[y_1,\ldots ,y_m,x_1,\ldots ,x_n]$. Consequently, basic
applications of Gr\"obner bases at the level of $A$ can be realized
by using Gr\"obner bases at the level of $T$, as presented in
Section 2.  Notations and conventions fixed in the previous sections
are maintained. {\parindent=0pt\v5

{\bf 4.1. Proposition} Let $R=E[a_1,\ldots ,a_m]$ be as given above.
If linear equations are solvable in $E$, then linear equations are
solvable in $R$. \vskip 6pt

{\bf Proof} If linear equations are solvable in $E$, then Gr\"obner
bases in the polynomial $E$-algebra $E[y_1,\ldots ,y_m]$ are
computable, thereby linear equations are solvable in $E[y_1,\ldots
,y_m]$ are solvable by ([AL], Chapter 4). Now, since $R\cong
E[y_1,\ldots ,y_m]/I$, the conclusion follows from Proposition
2.1.\QED}\v5

By Proposition 4.1 and Theorem 2.3, in principle  Gr\"obner bases in
the polynomial $R$-algebra $A=R[x_1,\ldots ,x_n]$ in any $n\ge 1$
variables are well defined, and Gr\"obner bases for ideal in $A$ can
be obtained by constructing Gr\"obner bases in the polynomial
$E[y_1,\ldots ,y_m]$-algebra $Q=(E[y_1,\ldots ,y_m])[x_1,\ldots
,x_n]$ in $n$-variables $x_1,\ldots ,x_n$ (i.e., by constructing
Gr\"obner bases with coefficients in $E[y_1,\ldots ,y_m]$). However,
with the aid of a reasult of [AB], we will show that Gr\"obner bases
for ideals in $A$ can indeed be obtained  by constructing Gr\"obner
bases in the polynomial $E$-algebra $T=E[y_1,\ldots ,y_m,x_1,\ldots
,x_n]$ in $m+n$ variables (i.e., by constructing Gr\"obner bases
with coefficients in $E$).
\par

To continue the discussion,  let us first bear in mind that as an
associative  ring,
$$T=E[y_1,\ldots ,y_m,x_1,\ldots ,x_n]=(E[y_1,\ldots ,y_m])[x_1,\ldots ,x_n]=Q.$$
For convenience, we write
$$\begin{array}{l}
\B =\{x^{\alpha}y^{\beta}=x_1^{\alpha_1}\cdots
x_n^{\alpha_n}y_1^{\beta_1}\cdots y_m^{\beta_m}~|~\alpha
=(\alpha_1,\ldots ,\alpha_n)\in\NZ^n,~\beta
=(\beta_1,\ldots ,\beta_m)\in\NZ^m\},\\
\B_x=\{x^{\alpha}=x_1^{\alpha_1}\cdots x_n^{\alpha_n}~|~\alpha
=(\alpha_1,\ldots ,\alpha_n)\in\NZ^n\},\\
\B_y=\{y^{\beta}=y_1^{\beta_1}\cdots y_m^{\beta_m}~|~\beta
=(\beta_1,\ldots ,\beta_m)\in\NZ^m\},\\
\end{array}$$ for the
standard free $E$-basis of $T$, the standard free $R$-basis of $A$,
and the standard free $E$-basis of $E[y_1,\ldots ,y_m]$,
respectively. Since any monomial ordering $\prec$ on $\B$ gives rise
to a monomial ordering on $\B_x$ and $\B_y$ respectively, we use the
same $\prec$ to denote the two obtained  monomial ordering on $\B_x$
and $\B_y$ respectively. Thus, if $0\ne f\in T$, then, as an element
of $Q$, $f$ can be expressed, with respect to $\prec$ on $\B_x$, as
$$f=h(y)x^{\alpha}+~\hbox{lower terms in the}~x~\hbox{variables},$$
where $h(y)\in E[y_1,\ldots ,y_m]$. So, as an element of $Q$, the
leading coefficient of $f$ with respect to $\prec$ on $\B_x$ is then
defined to be $\LC_x(f)=h(y)$.\v5

An indispensable bridge in reaching our main result of this section
(Theorem 4.4.(ii)) is the following  {\parindent=0pt\v5

{\bf 4.2. Proposition} ([AB], Proposition 2.1) With notation as
above, let $J$ be an ideal of the polynomial ring $T=E[y_1,\ldots
,y_m,x_1,\ldots ,x_n]=(E[y_1,\ldots ,y_m])[x_1,\ldots ,x_n]=Q$, and
let $\prec$ be a monomial ordering on $\B$ such that
$$x^{\alpha}y^{\beta}\prec x^{\alpha '}y^{\beta '}\Leftrightarrow x^{\alpha}\prec x^{\alpha '}~\hbox{or}~
x^{\alpha}=x^{\alpha '}~\hbox{and}~y^{\beta}\prec y^{\beta '},$$
i.e., $\prec$ is an elimination ordering with the $x$ variables
larger than the $y$ variables. The following statements hold.\par

(i) If, with respect to $\prec$ on $\B$, $\G =\{ g_1,\ldots g_t\}$
is a Gr\"obner basis of $J$ in $T$ (i.e., a Gr\"obner basis with
coefficients in $E$), then $\G$ is a Gr\"obner basis of $J$ in $Q$
with respect to $\prec$ on $\B_x$ (i.e., a Gr\"obner basis with
coefficients in $E[y_1,\ldots ,y_m]$). \par

(ii) With $\G$ as presented in (i), if we put $G=\{
h_i=\LC_x(g_i)~|~g_i\in\G\}$, then $G$ is a Gr\"obner basis for the
ideal $J_y=\langle\LC_x(f)~|~f\in J\rangle$ of $E[y_1,\ldots ,y_m]$
with respect to $\prec$ on $\B_y$.\par\QED\v5

{\bf Remark} The assertion (ii) of Proposition 4.2 will not be used
in the remaining part of this paper. }\v5

To better relate the main result of this section to Theorem 2.3, in
what follows we let $R= E[y_1,\ldots y_m]/I$ with  $I=\langle
\nu_1,\ldots ,\nu_s\rangle$, and we let $\varphi$: $Q=(E[y_1,\ldots
,y_m])[x_1,\ldots ,x_n]\r R[x_1,\ldots ,x_n]=A$ be the canonical
ring epimorphism with $\varphi (\sum_ih_ix^{\alpha (i)})=\sum_i\OV
h_ix^{\alpha (i)}$, where $h_i\in E[y_1,\ldots ,x_n]$, $x^{\alpha
(i)}\in\B_x$, and $\OV h_i=h_i+I$ is the coset in $R$ represented by
$h_i$. If $\prec$ is a given monomial ordering on $\B$, then, for
each $0\ne f\in T= Q$, we write $\LT (f)$ for the leading term of
$f$ in $T$ with respect to $\prec$ on $\B$, and we write $\LT_x(f)$
for the leading term of $f$ in $Q$ with respect to $\prec$ on
$\B_x$. Moreover, if $\varphi (f)=\OV f=\sum_i\OV h_ix^{\alpha
(i)}\ne 0$, we write $\LT_x(\OV f)$ for the leading term of $\OV f$
in $A$ with respect to $\prec$ on $\B_x$; and if $0\ne h\in
E[y_1,\ldots ,y_m]$, we write $\LT_y(h)$ for the leading term of $h$
with respect to $\prec$ on $\B_y$. Consequently we have the
corresponding $\LM (f)$, $\LM_x(f)$, $\LM_x(\OV f )$ and $\LM_y(h)$
in $T$, $Q$, $A$ and $E[y_1,\ldots ,y_m]$ respectively.
{\parindent=0pt\v5

{\bf 4.3. Lemma} Let $\prec$ be the elimination ordering on $\B$ as
described in Proposition 4.2. With notation as above, let $0\ne
f=\sum_{i,j}\varepsilon_{ij}x^{\alpha (i)}y^{\beta (j)}\in T$ with
$\LT (f)=\varepsilon_{k\ell}x^{\alpha (k)}y^{\beta (\ell )}$, where
$\varepsilon_{ij}\in E$ and $x^{\alpha (i)}y^{\beta (j)}\in\B$.  The
following statements hold.\par

(i) Viewing $f$ as an element of $Q$, if, after rewriting  $f$ in
the $x$ variables, $f=\sum_qh_qx^{\alpha (q)}$ with $h_q\in
E[y_1,\ldots ,y_m]$, then $\LT_x(f)=h_kx^{\alpha (k)}$ with
$\LT_y(h_k)=\varepsilon_{k\ell}y^{\beta (\ell )}$.\par

(ii) Let $h_k$ be as in (i). Then $f$ can be reduced modulo $\{
\nu_1,\ldots ,\nu_s\}$ in $T$ with respect to $\prec$ on $\B$ if and
only if $h_k$ can be reduced modulo $\{ \nu_1,\ldots ,\nu_s\}$ in
$E[y_1,\ldots ,y_m]$ with respect to $\prec$ on $\B_y$.\vskip 6pt

{\bf Proof} (i) Since as an element of $T$,
$f=\sum_{i,j}\varepsilon_{ij}x^{\alpha (i)}y^{\beta (j)}\in T$ with
$\LT (f)=\varepsilon_{k\ell}x^{\alpha (k)}y^{\beta (\ell )}$ by the
assumption, it follows from the expression $f=\sum_qh_qx^{\alpha
(q)}$ of $f$ in $Q$, and the definition of $\prec$ that
$\LT_x(f)=h_kx^{\alpha (k)}$ with
$\LT_y(h_k)=\varepsilon_{k\ell}y^{\beta (\ell )}$.}\par

(ii) Suppose that $f$ can be reduced modulo $\{ \nu_1,\ldots
,\nu_s\}$ in $T$ with respect to $\prec$ on $\B$. Then there exist
$\varepsilon_1,\ldots \varepsilon_s\in E$ and $x^{\alpha
(1)}y^{\beta (1)},\ldots ,x^{\alpha (s)}y^{\beta (s)}\in\B$ such
that
$$\begin{array}{l}
\LT (f)=\varepsilon_{k\ell}x^{\alpha (k)}y^{\beta (\ell
)}=\varepsilon_1x^{\alpha (1)}y^{\beta (1)}\LT (\nu_1)+\cdots +
\varepsilon_sx^{\alpha (s)}y^{\beta (s)}\LT
(\nu_s),\\
\LM (f)=x^{\alpha (k)}y^{\beta (\ell )}= x^{\alpha (p)}y^{\beta
(p)}\LM (\nu_p)~\hbox{for all}~\varepsilon_p\ne 0.
\end{array}\eqno{(1)}$$
Noticing that $\{\nu_1,\ldots ,\nu_s\}\subset I\subset E[y_1,\ldots
,y_m]$, the feature of $\prec$ then entails that $\LM
(\nu_p)=\LM_y(\nu_p)$ and $\LT (\nu_p)=\LT_y(\nu_p)$ for $1\le p\le
s$. It follows from (1) and (i) that
$$\begin{array}{l} x^{\alpha (k)}=x^{\alpha (p)}~\hbox{for all}~\varepsilon_p\ne 0,\\
\LM _y(h_k)=y^{\beta (\ell )}=y^{\beta (p)}\LM_y(\nu_p)~\hbox{for
all}~\varepsilon_p\ne 0,\\
\LT_y(h_k)=\varepsilon_{k\ell}y^{\beta (\ell
)}=\varepsilon_1y^{\beta (1)}\LT_y(\nu_1)+\cdots
+\varepsilon_sy^{\beta (s)}\LT_y(\nu_s).
 \end{array}\eqno{(2)}$$ This shows that
$h_k$ can be reduced modulo $\{\nu_1,\ldots ,\nu_s\}$ in
$E[y_1,\ldots ,y_m]$ with respect to $\prec$ on $\B_y$.\par

Conversely, suppose that $h_k$ can be reduced modulo $\{\nu_1,\ldots
,\nu_s\}$ in in $E[y_1,\ldots ,y_m]$ with respect to $\prec$ on
$\B_y$. Then there exist $\varepsilon_1,\ldots ,\varepsilon_s\in E$
and $y^{\beta (1)},\ldots ,y^{\beta (s)}\in\B_y$ such that
$$\begin{array}{l}
\LT_y(h_k)=\varepsilon_{k\ell}y^{\beta (\ell
)}=\varepsilon_1y^{\beta (1)}\LT_y(\nu_1)+\cdots
+\varepsilon_sy^{\beta (s)}\LT_y(\nu_s),\\
\LM _y(h_k)=y^{\beta (\ell )}=y^{\beta (p)}\LM_y(\nu_p)~\hbox{for
all}~\varepsilon_p\ne 0.\end{array}\eqno{(3)}$$
Multiplying both sides of the two equalities of (3) by $x^{\alpha
(k)}$, and noticing again the feature of $\prec$, we have
$$\begin{array}{l} \LT (f)=\varepsilon_{k\ell}x^{\alpha (k)}y^{\beta (\ell )}=\varepsilon_1x^{\alpha (k)}y^{\beta (1)}\LT (\nu_1)+\cdots +
\varepsilon_sx^{\alpha (k)}y^{\beta (s)}\LT (\nu_s),\\
\LM (f)=x^{\alpha (k)}y^{\beta (\ell )}=x^{\alpha (k)}y^{\beta
(p)}\LM (\nu_p)~\hbox{for all}~\varepsilon_p\ne 0.\end{array}$$
This shows that $f$ can be reduced modulo $\{ \nu_1,\ldots ,\nu_s\}$
in $T$ with respect to $\prec$ on $\B$, as desired.\par\QED
{\parindent=0pt\v5

{\bf 4.4. Theorem} With the ideal $I=\langle\nu_1,\ldots
,\nu_s\rangle\subset E[y_1,\ldots ,y_m]$ and $R=E[y_1,\ldots
,y_m]/I$ as before, given an ideal $\OV J=\langle \OV S\rangle$ of
$A=R[x_1,\ldots ,x_n]$ generated by $\OV S=\{ \OV f_1,\ldots ,\OV
f_t\}$, where $\OV f_i=\sum_j\OV h_{ij}x^{\alpha (j)}$ with $\OV
h_{ij}=h_{ij}+I$ the coset in $R$ represented by $h_{ij}\in
E[y_1,\ldots ,y_m]$, $1\le i\le t$, let $J=\langle S\rangle$ be the
ideal of $T=E[y_1,\ldots ,y_m,x_1,\ldots ,x_n]$ generated by $S=\{
f_1,\ldots ,f_t, \nu_1,\ldots ,\nu_s\}$, where
$f_i=\sum_jh_{ij}x^{\alpha (j)}$, $1\le i\le t$. Let $\prec$ be the
elimination  ordering on the standard free $E$-basis $\B$ of $T$ as
described in Proposition 4.2. Suppose that $G=\{ \nu_1,\ldots
,\nu_s\}$ forms a Gr\"obner basis of $I$ with respect to $\prec$ on
the standard free $E$-basis $\B_y$ of $E[y_1,\ldots ,y_m]$. The
following statements hold.\par

(i) If,  with respect to $\prec$ on $B$, $\G =\{ g_1,\ldots ,g_m\}$
is a Gr\"obner basis of $J$ constructed by means of ([AL], Algorithm
4.2.1) using the initial input-data $S$ in $T$, then for each
$g_i\not\in G$, we have $\LT_x (\varphi (g_i))=\varphi (\LT_x
(g_i))$ and hence $\LM_x (\varphi (g_i))=\LM_x (g_i)$, where
$$\varphi: \quad Q=(E[y_1,\ldots ,y_m])[x_1,\ldots
,x_n]\mapright{}{} R[x_1,\ldots ,x_n]=A$$ is the canonical ring
epimorphism with $\varphi (\sum_ih_ix^{\alpha (i)})=\sum_i\OV
h_ix^{\alpha (i)}$.\par

(ii) Let $\G$ be the Gr\"obner basis of $J$ presented in (i).
Considering the image $\varphi (\G )$ in $A$, the set
$$\OV{\G}=\{ \OV g_k=\varphi (g_k)~|~g_k\in\G ,~g_k\ne \nu_j,~1\le j\le s\}=\varphi (\G )-\{ 0\}$$
is a Gr\"obner basis for $\OV J$ with respect to the monomial
ordering $\prec$ on $\B_x$. \par

(iii) With $\G$ as presented in (i), each $g_k\in\G-\{\nu_1,\ldots
,\nu_s\}$, as an element of $Q$, has a representation
$g_k=g_k'+g_k''$ in which
$$\begin{array}{l} g_k'=\sum_qh_{kq}x^{\alpha (q)}~\hbox{with}~h_{kq}\not\in I~\hbox{for all}~q,\\
g_k''=\sum_{\ell}h_{k\ell}x^{\alpha (\ell
)}~\hbox{with}~h_{k\ell}\in I~\hbox{for all}~\ell,\end{array}$$ and
such a representation can be algorithmically determined. Hence each
element of the Gr\"obner basis $\OV{\G}$ obtained in (ii) has a
``real representation" in $A$, i.e., $\OV g_k=\varphi (g_k)=\varphi
(g_k')=\sum_q\OV{h}_{kq}x^{\alpha (q)}$ with all the $\OV{h}_{kq}\ne
0$. \vskip 6pt

{\bf Proof} (i) Let $g_i\in\G$ with $g_i\not\in G=\{ \nu_1,\ldots
,\nu_s\}$. If $g_i\in\{ f_1,\ldots ,f_t\}$, then it is clear that
that $\LT_x(\varphi (g_i))=\varphi (\LT_x(g_i))$. If $g_i\not\in\{
f_1,\ldots ,f_t,\nu_1,\ldots ,\nu_s\}$, then it follows from ([AL],
Algorithm 4.2.1) that $g_i$ is obtained by passing through some
$i$-th round executing the While loop, that is, $g_i:=r$, where $r$
appears as the remainder of a reduction
$$H_1g_1+\cdots +H_tg_t~\mapright{G'}{}_+r$$
which is minimal with respect to $G'$. So, $g_i$ cannot be further
reduced modulo $G'$. Note that $G\subset S\subset G'$ since $S$ is
the initial input-data running the algorithm.  Hence, $g_i$  cannot
be reduced modulo $G$ in $T$. It follows from Lemma 4.3 that if, as
an element of $T$, $g_i=\sum_{p,q}\varepsilon_{pq}x^{\alpha
(p)}y^{\beta (q)}$ with $\LT (g_i)=\varepsilon_{p'q'}x^{\alpha
(p')}y^{\beta (q')}$, and if, as an element of $Q$ expressed in the
$x$ variables, $g_i=\sum_kh_kx^{\alpha (k)}$ with $h_k\in
E[y_1,\ldots ,y_m]$, then $\LT_x(g_i)=h_{p'}x^{\alpha (p')}$ and
$h_{p'}$ cannot be reduced modulo $G$. Thereby we conclude that
$\varphi (g_i)=\OV h_{p'}x^{\alpha (p')}+\sum_{k\ne p'}\OV
h_kx^{\alpha (k)}$ with $\OV h_{p'}\ne 0$. Otherwise, $\OV h_{p'}=0$
would imply $h_{p'}\in I$ and hence $h_{p'}$ could be reduced modulo
$G$ (note that $G$ is a Gr\"obner basis of $I$). Therefore, we have
$\LT_x(\varphi (g_i))=\OV h_{p'}x^{\alpha (p')}=\varphi
(\LT_x(g_i))$, as desired. }\par

(ii) Since (i) holds and, by Proposition 4.2(i), $\G$ is a Gr\"obner
basis of $J$ in $Q$ with respect to $\prec$ on $\B_x$,  we are now
in the situation of Theorem 2.3 with $\Lambda =E[y_1,\ldots ,y_m]$.
So the proof of Theorem 2.3(ii) can be completely adapted to show
that $\OV{\G}$ is a Gr\"obner basis for $\OV J$ with respect to the
monomial ordering $\prec$ on $\B_x$.\par

(iii) Since $\G$ is finite and since $G=\{\nu_1,\ldots,\nu_s\}$ is a
Gr\"obner basis of $I$ by our assumption, it follows that there is
an algorithm to determine whether a coefficient of $g_k$ (as an
element of $Q$) is in $I=\langle\nu_1,\ldots ,\nu_s\rangle$. \QED\v5

In light of Proposition 4.2, one may check that the following
holds.{\parindent=0pt\v5

{\bf 4.5. Theorem} Let the Gr\"obner bases $\G$ and $\OV{\G}$ be as
in Theorem 4.4.  Then, similar results as presented in Proposition
2.4 -- Theorem 2.8 hold, that is, all basic applications of the
Gr\"obner basis $\OV{\G}$ at the level of $A=(E[y_1,\ldots
,y_m]/I)[x_1,\ldots ,x_n]$ can be realized by using the Gr\"obner
basis $\G$ at the level of $T=E[y_1,\ldots ,y_m,x_1,\ldots ,x_n]$.}

\section*{5. The Case of $R =D[\vartheta_1,\ldots ,\vartheta_m]$ with $D$ a PID}
In consideration of the prominent feature of strong Gr\"obner bases
over a PID (see [AL], Section 4.5), in this section we demonstrate
in more details the application of Theorem 4.4 and Theorem 4.5 to a
finitely generated $D$-algebra $R =D[\vartheta_1,\ldots
,\vartheta_m]$ with the generating set $\{ \vartheta_1,\ldots
,\vartheta_m\}$, where $D$ is a PID including the case of $D=K$
being a field. Since $R\cong D[y_1,\ldots ,y_m]/I$, where
$D[y_1,\ldots ,y_m]$ is the polynomial ring in $m$ variables over
$D$ and $I$ is an ideal of $D[y_1,\ldots ,y_m]$, in what follows we
let $R=D[y_1,\ldots ,y_m]/I$. All notations and conventions used
before are maintained.{\parindent=0pt\v5


{\bf 5.1. Corollary} If  in Theorem 4.4 and Theorem 4.5, $E=D$ is a
PID, respectively $E=K$ is a field, then
\par

(i) the Gr\"obner basis $\OV{\G}$ in the sense of Definition 1.3 can
be obtained by constructing a strong Gr\"obner basis $\G$ in the
polynomial ring $D[y_1,\ldots ,y_m,x_n,\ldots ,x_n]$ by means of
([AL], Algorithm 4.5.1), respectively by constructing a Gr\"obner
basis $\G$ in the polynomial ring $K[y_1,\ldots ,y_m,x_n,\ldots
,x_n]$ by means of the classical Buchberger algorithm; and\par

(ii) similar results as presented in Proposition 2.4 -- Theorem 2.8
hold, that is, all basic applications of the Gr\"obner basis
$\OV{\G}$ at the level of $A=(D[y_1,\ldots ,y_m]/I)[x_1,\ldots
,x_n]$ can be realized by using the strong Gr\"obner basis $\G$ at
the level of $T=D[y_1,\ldots ,y_m,x_1,\ldots ,x_n]$, respectively
all basic applications of the Gr\"obner basis $\OV{\G}$ at the level
of $A=(K[y_1,\ldots ,y_m]/I)[x_1,\ldots ,x_n]$ can be realized by
using the (classical)  Gr\"obner basis $\G$ at the level of
$T=K[y_1,\ldots ,y_m,x_1,\ldots ,x_n]$.\par\QED}\v5

In the case of $E=K$ being a field, the following two significant
cases immediately illustrate the advantage of Corollary 5.1:\par

(1) $R=K[y_1,\ldots ,y_m]/M$ is an extension field of $K$, where
$M=\langle f_1,\ldots f_t\rangle$ is a maximal ideal of
$K[y_1,\ldots ,y_m]$.\par

(2) $R=K[f_1,\ldots ,f_m]$ is the $K$-subalgebra generated by
polynomials  $f_1,\ldots ,f_m\in K[x_1,\ldots ,x_n]$. Note that  in
this case $R\cong K[y_1,\ldots ,y_m]/I$, and that by ([AL], Theorem
2.4.2), a Gr\"obner basis of $I$ can be worked out by calculating a
Gr\"obner basis of $I=H\cap K[y_1,\ldots ,y_m]$ by means of the
classical Buchberger algorithm, where $H=\langle y_1-f_1,\ldots
,y_m-f_m\rangle\subset K[y_1,\ldots ,y_m,x_1,\ldots
,x_n]$.{\parindent=0pt\v5

{\bf Remark} Note that every Gr\"obner basis over a field $K$ is
certainly a strong Gr\"obner basis in the sense of Definition 3.1,
and note also that in Corollary 5.1(i) we emphasized that the
obtained Gr\"obner $\OV{\G}$ is the one in the sense of Definition
1.3 instead of a strong Gr\"obner basis. The reason is that if $m\ge
2$, then, as it is pointed out to us by [AL] on page 251, the
Gr\"obner basis $\G$ constructed in $T=K[y_1,\ldots ,y_m,x_1,\ldots
,x_n]$ may not be a strong Gr\"obner basis in $Q=(K[y_1,\ldots
y_m])[x_1,\ldots ,x_n]$. So, the Gr\"obner basis $\OV{\G}$ we
obtained in the way as presented in Theorem 4.4 may not be a strong
Gr\"obner basis. Nevertheless, this does not matter practical
applications of Gr\"obner bases in $A$, because Theorem 5.1(ii)
tells us that the basic applications of $\OV{\G}$ are realized by
using the (strong or classical) Gr\"obner bases $\G$ at the level of
$T$. }\par

Whereas we will see that in the case that $K[y]$ is the polynomial
ring in one variable $y$ over a field $K$, Corollary 5.1 may be
turned to be much better due to a nice result of [AL] quoted below.
{\parindent=0pt\v5

{\bf 5.2. Proposition} ([AL], Theorem 4.5.12) Let $K[y]$ be the
polynomial ring in one variable $y$ over a field $K$. Then $\G =\{
g_1,\ldots ,g_m\}$ is a Gr\"obner basis in $T=K[y,x_1,\ldots ,x_n]$
with respect to an elimination ordering with the $x$ variables
larger than $y$ if and only if $\G$ is a strong Gr\"obner basis in
$Q=(K[y])[x_1,\ldots ,x_n]$.\par\QED\v5

{\bf 5.3. Theorem} With notation as fixed above, let $I=\langle
\nu\rangle$ be an ideal of $K[y]$ generated by the single polynomial
$\nu$, and $R=K[y] /I$.  Let $\OV J=\langle\OV S\rangle$ be an ideal
of $A=R[x_1,\ldots ,x_n]$ generated by the set of nonzero elements
$\OV S=\{\OV f_1,\ldots ,\OV f_s\}$, where $\OV
f_i=\sum_j\OV{\lambda}_{ij}x^{\alpha (j)}$, $1\le i\le s$. Consider
in $T=K[y,x_1,\ldots ,x_n]$ the set of elements $S=\{ f_1,\ldots
,f_s,\nu\}$, where $f_i=\sum_j\lambda_{ij}x^{\alpha (j)}$, and let
$J=\langle S\rangle$ be the ideal of $T$ generated by $S$.  If, with
respect to an elimination ordering $\prec$ with $y$ smaller than all
the $x_i$, $\G =\{ g_1,\ldots ,g_m\}$ is a  Gr\"obner basis of $J$
constructed by the classical Buchberger algorithm using the initial
input-data $S$ in $T$, then  \par

(i) the set
$$\OV{\G}=\{ \OV g_k=\varphi (g_k)~|~g_k\in\G ,~g_k\ne a\}=\varphi (\G )-\{ 0\}$$
is a strong Gr\"obner basis for $\OV J$ in the sense of Definition
3.1, where $\varphi$ is the canonical ring epimorphism from
$Q=(K[y])[x_1,\ldots ,x_n]$ to $A$; \par

(ii) each $g_k\in\G-\{\nu\}$, as an element of $Q$, has a
representation $g_k=g_k'+g_k''$ in which
$$\begin{array}{l} g_k'=\sum_qh_{kq}x^{\alpha (q)}~\hbox{with}~h_{kq}\not\in I~\hbox{for all}~q,\\
g_k''=\sum_{\ell}h_{k\ell}x^{\alpha (\ell
)}~\hbox{with}~h_{k\ell}\in I~\hbox{for all}~\ell,\end{array}$$ and
such a representation can be algorithmically determined, thereby
each element of the Gr\"obner basis $\OV{\G}$ obtained in (ii) has a
``real representation" in $A$, i.e., $\OV g_k=\varphi (g_k)=\varphi
(g_k')=\sum_q\OV{h}_{kq}x^{\alpha (q)}$ with all the $\OV{h}_{kq}\ne
0$; and\par

(iii) similar results as presented in Proposition 2.4 -- Theorem 2.8
hold, that is, all basic applications of the Gr\"obner basis
$\OV{\G}$ at the level of $A$ can be realized by using the
(classical) Gr\"obner basis $\G$ at the level of $T$.\vskip 6pt

{\bf Proof} (i) By Proposition 5.2, $\G$ is first of all a strong
Gr\"obner basis of $J$ in $Q$. Now, the remaining argument is
similar to the proof of Theorem 3.2.\par

(ii) Since $\G$ is finite and since $I=\langle\nu\rangle\subset
K[y]$, where $K[y]$ is the polynomial ring in one variable over a
field $K$, the classical division algorithm can be used to determine
whether a coefficient of $g_k$ (as an element of $Q$) is in $I$.\par

(iii) This follows from Corollary 5.1(ii).\QED}\v5

We next consider the case when $D$ is a PID with the field of
fractions $K$. Let $K\subset L$ be a field extension and let
$\vartheta\in L$ be an algebraic element over $K$. We are interested
in the ring extension $D\subset R=D[\vartheta ]$ where the minimal
polynomial of $\vartheta$ over $K$ is known. {\parindent=0pt\v5

{\bf 5.4. Lemma} Let $R=D[\vartheta ]$ be as above. Suppose that
$q(y)\in K[y]$ is the minimal polynomial of $\vartheta$ over $K$.
Then there are $b,d\in D$ and a primitive polynomial $p(y)\in D[y]$
(which is necessarily irreducible) such that $bq(y)=dp(y)$. Moreover
$D[\vartheta ]\cong D[y]/\langle p(y)\rangle$.\par\QED\v5

{\bf 5.5. Theorem} With $R=D[\vartheta ]$ as fixed above, let $\OV
J=\langle\OV S\rangle$ be an ideal of $A=R[x_1,\ldots ,x_n]$
generated by the set of nonzero elements $\OV S=\{\OV f_1,\ldots
,\OV f_s\}$, where $\OV f_i=\sum_j\lambda_{ij}(\vartheta )x^{\alpha
(j)}$ with $\lambda_{ij}(\vartheta )\in R$, $1\le i\le s$. Let
$p(y)$ be the polynomial as presented in Lemma 5.4. Consider in the
polynomial ring $T=D[y, x_1,\ldots ,x_n]$ the set of elements $S=\{
f_1,\ldots ,f_s,p(y)\}$, where $f_i=\sum_j\lambda_{ij}(y)x^{\alpha
(j)}$ with $\lambda_{ij}(y)$ given by replacing $\vartheta$ by $y$
in the coefficients of $\OV f_i$, and let $J=\langle S\rangle$ be
the ideal of $T$ generated by $S$.  If, with respect to a monomial
ordering $\prec$ as described in Proposition 4.2, $\G =\{ g_1,\ldots
,g_m\}$ is a strong Gr\"obner basis of $J$ constructed by means of
([AL], Algorithm 4.5.1) using the initial input-data $S$ in $T$,
then\par

(i)  the set
$$\OV{\G}=\{ \OV g_k=\varphi (g_k)~|~g_k\in\G ,~g_k\ne a\}=\varphi (\G )-\{ 0\}$$
is a Gr\"obner basis for $\OV J$ in the sense of Definition 1.3,
where $\varphi$: $T\r A$ is the canonical ring epimorphism with
$\varphi (\sum_J\lambda_j(y)x^{\alpha
(j)})=\sum_j\lambda_j(\vartheta)x^{\alpha (j)}$;\par

(ii) each $g_k\in\G-\{p(y)\}$, as an element of $Q$, has a
representation $g_k=g_k'+g_k''$ in which
$$\begin{array}{l} g_k'=\sum_q\lambda_{kq}(y)x^{\alpha (q)}~\hbox{with}~\lambda_{kq}(y)\not\in \langle p(y)\rangle~\hbox{for all}~q,\\
g_k''=\sum_{\ell}\lambda_{k\ell}(y)x^{\alpha (\ell
)}~\hbox{with}~\lambda_{k\ell}(y)\in \langle p(y)\rangle~\hbox{for
all}~\ell,\end{array}$$ and such a representation can be
algorithmically determined, thereby  each element of the Gr\"obner
basis $\OV{\G}$ obtained in (ii) has a ``real representation" in
$A$, i.e., $\OV g_k=\varphi (g_k)=\varphi
(g_k')=\sum_q\lambda_{kq}(\vartheta )x^{\alpha (q)}$ with all the
$\lambda_{kq}(\vartheta )\ne 0$; and\par

(iii) similar results as presented in Proposition 2.4 -- Theorem 2.8
hold, that is, all basic applications of the Gr\"obner basis
$\OV{\G}$ at the level of $A$ can be realized by using the strong
Gr\"obner basis $\G$ at the level of $T$.\vskip 6pt

{\bf Proof} (i) Since a strong Gr\"obner basis in the sense of
Definition 3.1 is certainly a Gr\"obner basis in the sense of
Definition 1.3, it follows from Theorem 4.4 that $\OV{\G }$ is a
Gr\"obner basis for the ideal $\OV J$.\par

(ii) Since $\G$ is finite and since $\langle p(y)\rangle\subset
D[y]$, where $D[y]$ is the polynomial ring in one variable over the
PID $D$ in which linear equitions are solvable, it follows that
there is an algorithm to determine whether a coefficient of $g_k$
(as an element of $Q$) is in $I=\langle\nu_1,\ldots
,\nu_s\rangle$.\par

(iii) This follows from Theorem 4.5. \QED\v5

{\bf Remark} Note that the Gr\"obner basis $\G$ in Theorem 5.5 is a
strong Gr\"obner basis in the sense of Definition 3.1, and note also
that we emphasized that the obtained Gr\"obner $\OV{\G}$ is the one
in the sense of Definition 1.3 instead of a strong Gr\"obner basis.
The reason is that the Gr\"obner basis $\G$  may not be a strong
Gr\"obner basis in $Q=(D[y])[x_1,\ldots ,x_n]$, even if
$D=\mathbb{Z}$ is the ring of integers ([AL], Example 4.5.4). So,
the Gr\"obner basis $\OV{\G}$ we obtained in the way as presented in
Theorem 5.5(i) may not be a strong Gr\"obner basis. Nevertheless,
this does not matter practical applications of Gr\"obner bases in
$A$, because Theorem 5.5.(ii) tells us that the basic applications
of $\OV{\G}$ are realized by using the strong Gr\"obner basis $\G$
at the level of $T$. }\v5

Obviously, Theorem 5.5 immediately applies to the case of $R =
\mathbb{Z}[\vartheta ]$, where $\vartheta$ is an {\it arbitrary
algebraic number}, thereby $\mathbb{Z}[\vartheta ]$ may not even be
a UFD.  \v5

\section*{6. The Case of $R=\mathbb{Z}_{p^n}[y]/\langle f\rangle$ Being a Galois Ring}
We end this paper by applying the results of previous sections to
the Gr\"obner basis theory over Galois rings.  \v5

Let $p$ be a prime in  $\mathbb{Z}$ and $m$, $n$ be positive
integers. Recall from the literature ([Mc], [Rag]) that if $\hat f$
is a monic basic irreducible polynomial of degree $m$ in the
one-variable polynomial ring $\mathbb{Z}_{p^n}[y]$ (i.e., $\hat f$
is irreducible modulo $p$), then   the quotient ring
$R=\mathbb{Z}_{p^n}[y]/\langle \hat f\rangle$ is called the Galois
ring of order $p^{mn}$ and characteristic $p^n$.\par

Consider the polynomial ring $A=R[x_1,\ldots ,x_n]$. For
convenience, if $\hat f=\sum_i\bar a_iy^i\in\mathbb{Z}_{p^n}[y]$,
where each $a_i\in\mathbb{N}$ satisfies $1\le a_i\le p^n-1$, then we
write $f=\sum_ia_iy^i$ for the corresponding polynomial in
$T=\mathbb{Z}[y,x_1,\ldots ,x_n]$; If $\OV J=\langle \OV S\rangle$
is an ideal of $A$ generated by the set of nonzero polynomials $S=\{
\OV f_1,\ldots ,\OV f_t\}$, where $\OV f_i=\sum_j\OV h_{ij}x^{\alpha
(j)}$ with $\OV h_{ij}\in R$ represented by $\hat h_{ij}=\sum_k\bar
a^{ij}_ky^k\in\mathbb{Z}_{p^n}[y]$ with $a^{ij}_k\in\NZ$ satisfying
$1\le a^{ij}_k\le p^n-1$, then we write $h_{ij}=\sum_ka^{ij}_ky^k$
and $f_i=\sum_jh_{ij}x^{\alpha (j)}$, $1\le i\le t$, and let
$J=\langle S\rangle$ be the ideal of $T$ generated by $S=\{
f_1,\ldots ,f_t,f,p^n\}$. {\parindent=0pt\v5

{\bf 6.1. Theorem} With notation as fixed above, if, with respect to
a monomial ordering $\prec$ as described in Proposition 4.2, $\G =\{
g_1,\ldots ,g_m\}$ is a strong Gr\"obner basis of $J$ constructed by
means of ([AL], Algorithm 4.5.1) using the initial input-data $S$ in
$T$, then \par

(i) the set
$$\OV{\G}=\{ \OV g_k=\varphi (g_k)~|~g_k\in\G ,~g_k\ne f\}=\varphi\psi (\G )-\{ 0\}$$
is a Gr\"obner basis for $\OV J$ in the sense of Definition 1.3,
where $\varphi$, and $\psi$ are the canonical ring epimorphisms as
shown below:
$$\mathbb{Z}[y,x_1,\ldots ,x_n]~\mapright{\psi}{}~\mathbb{Z}_{p^n}[y,x_1,\ldots ,x_n]
=(\mathbb{Z}_{p^n}[y])[x_1,\ldots
,x_n]~\mapright{\varphi}{}~(\mathbb{Z}_{p^n}[y]/\langle\hat
f\rangle)[x_1,\ldots ,x_n];$$\par

(ii) each $g_k\in\G-\{\nu\}$, as an element of $Q$, has a
representation $g_k=g_k'+g_k''$ in which
$$\begin{array}{l} g_k'=\sum_qh_{kq}x^{\alpha (q)}~\hbox{with}~h_{kq}\not\in \langle\hat
f\rangle~\hbox{for all}~q,\\
g_k''=\sum_{\ell}h_{k\ell}x^{\alpha (\ell
)}~\hbox{with}~h_{k\ell}\in \langle\hat f\rangle~\hbox{for
all}~\ell,\end{array}$$ and such a representation can be
algorithmically determined, thereby  each element of the Gr\"obner
basis $\OV{\G}$ obtained in (ii) has a ``real representation" in
$A$, i.e., $\OV g_k=\varphi (g_k)=\varphi
(g_k')=\sum_q\OV{h}_{kq}x^{\alpha (q)}$ with all the $\OV{h}_{kq}\ne
0$; and\par

(iii) similar results as presented in Proposition 2.4 -- Theorem 2.8
hold, that is, all basic applications of the Gr\"obner basis
$\OV{\G}$ at the level of $A$ can be realized by using the strong
Gr\"obner basis $\G$ at the level of $T$.

{\bf Proof} This follows from Theorem 3.2, Theorem 4.4 and Theorem
4.5.\QED\v5

{\bf Remark} (i)  All results obtained in this paper for ideals may
be generalized to modules without much difficulty.\par

(ii) In a forthcoming paper [Li], the methods proposed in this paper 
will be generalized to deal with Gr\"obner bases in certain 
noncommutative algebras over rings, in particular, the Gr\"obner 
bases in a solvable polynomial algebra (in the sense of [K-RW]) over 
a commutative ring, and the Gr\"obner bases in a free algebra 
$R\langle X_1,\ldots ,X_n\rangle$ over a quotient ring 
$R=K[y_1,\ldots ,y_m]/I$ of a commutative polynomial ring 
$K[y_1,\ldots ,y_m]$ (including $I=\{ 0\}$).

\parindent=0pt\vskip 1truecm

\centerline{References}\parindent=1truecm\par

\item{[AB]} W.W. Adams and A. K. Boyle, Some results on Gr\"obner bases over commutative rings,
{\it J. Symbolic Comp.}, 13(1992), 473 -- 484.

\item{[AL]} W.W. Adams and P. Loustaunau, {\it An Introduction to Gr\"obner Bases},
Graduate Studies in Mathematics, Vol. 3. American Mathematical
Society, 1994.

\item{[BF]} E. Byrne and P. Fitzpatrick, Gr\"obner bases over Galois rings with an application
to decoding alternant codes, {\it J. Symbolic Comp.}, 31(2001), 565
-- 584.

\item{[BM]} E. Byrne and T. Mora, Gr\"obner Bases over Commutative Rings and Applications to Coding
Theory Gr?bner Bases,  In: Sala, Mora, Perret, Sakata, Traverso
(eds)., {\it Groebner Bases, Coding, and Cryptography},  RISC
Series: Springer (2009), 239-262.

\item{[Bu1]} B. Buchberger, {\it Ein Algorithmus zum Auffinden der
Basiselemente des Restklassenringes nach einem nulldimensionalen
polynomideal}, PhD thesis. University of Innsbruck, 1965.

\item{[Bu2]} B. Buchberger, Gr\"obner bases: An algorithmic method in
polynomial ideal theory, in: {\it Multidimensional Systems Theory}
(N.K. Bose, ed.). Reidel Dordrecht, 1985, 184--232.

\item{[BW]} T. Becker and V. Weispfenning, {\it Grobner bases. A computational approach to commutative
algebra}, Graduate Texts in Mathematics 141. Springer-Verlag, New
York, 1993.

\item{[GTZ]} P. Gianni, B. Trager, G. Zacharias, Gr¡§obner bases and primary
decomposition of polynomial ideals, {\it J. Symbolic Comp.}, 6(
1988), 148-166, a.k.a. Computational Aspects of Commutative Algebra
(L. Robbiano, ed.). Academic Press, San Diego, 1989, 15-33. ISBN
0-12-589590-9.

\item{[KC]} D. Kapur and Y. Cai, An algorithm for computing a Gr\"obner basis of a polynomial ideal
over a ring with zero divisors, {\it Preprint}, University of New
Mexico, 2003. http://www.cs.unm.edu/~treport/tr/03-12/GB.pdf

\item{[K-RK]} A. Kandri-Rody and D. Kapur , Computing the Gr\"obner
basis of an ideal in polynomial rings over a Euclidean ring, {\it J.
Symbolic Comp.}, 6(1990), 37 -- 56.

\item{[K-RW]} A.~Kandri-Rody and V.~Weispfenning, Non-commutative
Gr\"obner bases in algebras of solvable type, {\it J. Symbolic
Comput.}, 9(1990), 1--26.

\item{[Li]} Huishi Li, On the construction of noncommutative Gr\"obner bases with coefficients in quotient rings, 
{\it in preparation}.

\item{[Mc]} B.R. McDonald, {\it Finite Rings with Identity}, New York,
Marcel Dekker, 1974.

\item{[M\"ol]} H.M. M\"oller, On the construction of Gr\"obner bases
using syzygies, {\it J. Symbolic Comp.}, 6(1988), 345 -- 359.

\item{[NS1]} G.H. Norton and A. Salagean, On the structure of linear and cyclic codes over finite
chain rings, {\it Appl. Algebra Engrg. Comm. Comput}. 10(2000),
489-506.

\item{[NS2]} G.H. Norton and A. Salagean, Strong Gr\"obner bases for polynomials over a
principal ideal ring, {\it Bull. Austral. Math. Soc.}, 64(2001) 505
-- 528.

\item{[Pan]} L. Pan, On the D-bases of polynomial ideals over principal ideal domains, {\it J. Symbolic Comp.},
 7(1988), 55 -- 69.

\item{[Pau]} F. Pauer,  Gr\"obner bases with coefficients in rings,
{\it J. Symbolic Computation}, 42(2007), 1003 -- 1011.

\item{[Rag]} R. Raghavendran, Finite associative rings, {\it Compositio
Mathematica}, 21(1969), 195 -- 229.

\item{[Seel]} F. Seelisch, The theory of Gr\"obner bases over
coefficient rings and its application in electrical engineering,
2010. Talk slides at \par
http://www.mathematik.uni-kl.de/~seelisch/Talks/S2AM2010Berlin.pdf

\item{[Zac]} G. Zacharias, {\it Generalized Gr\"obner bases in commutative
polynomial rings}, Bachelor¡¯s thesis, M.I.T. 1978.

\end{document}